\documentclass{article}%
\usepackage{bbm}
\usepackage{epsfig}
\usepackage{graphics,graphicx,amssymb,amsmath,verbatim}
\usepackage{amssymb}
\usepackage{mathrsfs}
\usepackage{amsfonts}
\usepackage{amsmath}
\usepackage{graphicx}
\usepackage{easybmat}
\usepackage{color}%
\providecommand{\U}[1]{\protect \rule{.1in}{.1in}}
\newtheorem{theorem}{Theorem}

\newtheorem{assumption}{Assumption}

\newtheorem{corollary}{Corollary}

\newtheorem{lemma}{Lemma}

\newtheorem{proposition}{Proposition}
\newtheorem{remark}{Remark}

\allowdisplaybreaks[4]
\begin{document}

\title{Solutions to Linear Bimatrix Equations with Applications to Pole Assignment of
Complex-Valued Linear Systems}
\author{Bin Zhou\thanks{Center for Control Theory and Guidance Technology, Harbin
Institute of Technology, Harbin, 150001, China. Email: binzhoulee@163.com,
binzhou@hit.edu.cn.}}
\date{}
\maketitle

\begin{abstract}
We study in this paper solutions to several kinds of linear bimatrix equations
arising from pole assignment and stability analysis of complex-valued linear
systems, which have several potential applications in control theory,
particularly, can be used to model second-order linear systems in a very dense
manner. These linear bimatrix equations include generalized Sylvester bimatrix
equations, Sylvester bimatrix equations, Stein bimatrix equations, and
Lyapunov bimatrix equations. Complete and explicit solutions are provided in
terms of the bimatrices that are coefficients of the equations/systems. The
obtained solutions are then used to solve the full state feedback pole
assignment problem for complex-valued linear system. For a special case of
complex-valued linear systems, the so-called antilinear system, the solutions
are also used to solve the so-called anti-preserving (the closed-loop system
is still an antilinear system) and normalization (the closed-loop system is a
normal linear system) problems. Second-order linear systems, particularly, the
spacecraft rendezvous control system, are used to demonstrate the obtained
theoretical results.

\vspace{0.3cm}


\end{abstract}

\section{Introduction}

In this paper we continue to study complex-valued linear systems introduced in
\cite{zhou17arxiv}. Complex-valued linear systems refer to linear systems
whose state evolution dependents on both the state and its conjugate (see
Subsection \ref{sec2.1} for a detailed introduction). There are several
reasons for study this class of linear systems \cite{zhou17arxiv}, for
example, they are naturally encountered in linear dynamical quantum systems
theory, and can be used to model any real-valued linear systems with lower
dimensions (see Subsection \ref{sec3.3} for a detailed development). Analysis
and design of complex-valued linear systems have been studied in our early
paper \cite{zhou17arxiv}, where some fundamental problems such as state
response, controllability, observability, stability, pole assignment, linear
quadratic regulation, and state observer design, were solved. The conditions
and/or methods obtained there are based on bimatrices associated with the
complex-valued linear system, which is mathematically appealing.

The pole assignment problem for complex-valued linear system was solved in
\cite{zhou17arxiv} by using coefficients of the so-called real-representation
system, for which any pole assignment algorithms for normal linear systems can
be applied. In this paper, we will continue to study the pole assignment
problem for complex-valued linear systems by establishing a different method.
Our new solution is based on solving the so-called (generalized) Sylvester
bimatrix equation whose coefficients are bimatrices associated with the
complex-valued linear system. Our study on linear bimatrix equations and their
applications in pole assignment has been inspired by early work for normal
linear systems. For example, the (normal) Sylvester matrix equation was
utilized to solve the pole assignment for normal linear system in
\cite{bs82scl}, and the generalized Sylvester matrix equation was used in
\cite{duan93tac,duan96tac,duan15book,dz06tac} and \cite{zd06scl} to solve the
(parametric) pole assignment problem for normal linear systems, descriptor
linear systems, and even high-order linear systems.

We will show that the pole assignment problem for a complex-valued linear
system has a solution if and only if the associated (generalized) Sylvester
bimatrix equation has a nonsingular solution. Thus the main task of this paper
is to provide complete and explicit solutions to the homogeneous (generalized)
Sylvester bimatrix equation. The solutions we provided have quite element
expressions that use the original coefficient bimatrices and a right-coprime
factorization (in the bimatrix framework) of the system. We also provide
solutions to non-homogeneous Sylvester bimatrix equations and Stein bimatrix
equations which include the Lyapunov bimatrix equation as a special case.

We are particularly interested in pole assignment for the so-called antilinear
system studied recently in \cite{wdl13aucc,wqls16jfi, wz17book} and
\cite{wzls15iet}. By our approach we first provide closed-form solutions to
the associated (generalized) Sylvester bimatrix equations, and then consider
two different problems, namely, the anti-preserving problem which ensures that
the closed-loop system is still (or equivalent to) an antilinear system, and
the normalization problem which guarantees that the closed-loop system is
(equivalent to) a normal linear system. The anti-preserving problem was
firstly studied in \cite{wz17book}. However, we can provide complete solutions
that use full state feedback rather than only normal state feedback used in
\cite{wz17book}. We discovered that the anti-preserving problem is meaningful
only for discrete-time antilinear systems (as studied in \cite{wz17book})
since any continuous-time antilinear system cannot be asymptotically stable.
However, the normalization problem is valid for both continuous-time and
discrete-time antilinear systems, and seems more interesting as the
closed-loop system is (equivalent to) a normal linear system that is more easy
to handle.

\textbf{Notation}: For a matrix $A\in \mathbf{C}^{n\times m}$, we use $A^{\#},$
$A^{\mathrm{T}},$ $A^{\mathrm{H}},$ $\mathrm{rank}\left(  A\right)  ,$
$\left \vert A\right \vert ,$ $\left \Vert A\right \Vert ,$ $\lambda \left(
A\right)  ,$ $\rho \left(  A\right)  ,$ $\mu \left(  A\right)  ,$
$\operatorname{Re}\left(  A\right)  $ and $\operatorname{Im}\left(  A\right)
$ to denote respectively its conjugate, transpose, conjugate transpose, rank,
determinant (when $n=m$), norm, eigenvalue set (when $n=m$), spectral radius
(when $n=m$), spectral abscissa ($\max_{s\in \lambda(A)}\{ \operatorname{Re}%
\{s\} \}$), real part and imaginary part. The notation $0_{n\times m}$ refers
to an $n\times m$ zero matrix. For two integers $p,q$ with $p\leq q,$ denote
$\mathbf{I}\left[  p,q\right]  =\{p,p+1,\cdots,q\}.$ Let $\mathbf{Z}%
^{+}=\{0,1,2,\cdots \}$, $\mathbf{R}^{+}=[0,\infty),\mathbf{R}=\mathbf{R}%
^{+}\cup \{-\mathbf{R}^{+}\},\mathbf{Z}=\mathbf{Z}^{+}\cup \{-\mathbf{Z}^{+}\}$,
and $\mathrm{j}$ the unitary imaginary number. For a series of matrices
$A_{i},i\in \mathbf{I}\left[  1,l\right]  ,$ $\mathrm{diag}\{A_{1},A_{2}%
,\cdots,A_{l}\}$ denotes a diagonal matrix whose diagonal elements are
$A_{i},i\in \mathbf{I}\left[  1,l\right]  .$ 

\section{\label{sec2}Motivation and Preliminaries}

\subsection{\label{sec2.1}Complex-Valued Linear Systems}

To introduce complex valued linear systems we recall the bimatrix $\left \{
A_{1},A_{2}\right \}  \in \{ \mathbf{C}^{n\times m},\mathbf{C}^{n\times m}\}$
given in \cite{zhou17arxiv}, where it was defined in such a manner that, for
any $x\in \mathbf{C}^{m},$
\[
y=\left \{  A_{1},A_{2}\right \}  x\triangleq A_{1}x+A_{2}^{\#}x^{\#},
\]
which defines a linear mapping over the field of real numbers
\cite{zhou17arxiv}. Further properties of the bimatrix can be found in \cite{zhou17arxiv}. With the notion of bimatrix, we
continue to study the following complex-valued linear system
\cite{zhou17arxiv} (without output equation)%
\begin{equation}
x^{+}=\left \{  A_{1},A_{2}\right \}  x+\left \{  B_{1},B_{2}\right \}  u,
\label{sys}%
\end{equation}
where $A_{i}\in \mathbf{C}^{n\times n},B_{i}\in \mathbf{C}^{n\times m},i=1,2,$
are known coefficients, $x=x(t)$ is the state, $u=u(t)$ is the control, and
$x^{+}(t)$ denotes $x\left(  t+1\right)  $ if $t\in \mathbf{Z}^{+}$ (namely,
discrete-time systems) and denotes $\dot{x}(t)$ if $t\in \mathbf{R}^{+}$
(namely, continuous-time systems). Throughout this paper, the dependence of
variables on $t$ will be suppressed unless necessary. The initial condition is
set to be $x\left(  0\right)  =x_{0}\in \mathbf{C}^{n}$ \cite{zhou17arxiv}.

There are several reasons for studying linear systems in the form of
(\ref{sys}). Readers are encouraged to refer to \cite{zhou17arxiv} for
details, while an explicit application of (\ref{sys}) to second-order linear
system will be shown in detail in the next subsection. If $A_{2}\ $and $B_{2}$
are null matrices, then system (\ref{sys}) becomes
\begin{equation}
x^{+}=A_{1}x+B_{1}u, \label{normal}%
\end{equation}
which is the normal linear system that has been well studied during the past
half century \cite{kfa69book, rugh96book}. If $A_{1}$ and $B_{1}$ are set as
zeros, then (\ref{sys}) reduces to the so-called antilinear system
\begin{equation}
x^{+}=A_{2}^{\#}x^{\#}+B_{2}^{\#}u^{\#}, \label{antilinear}%
\end{equation}
which was initially studied in \cite{wdl13aucc} and \cite{wzls15iet}.

In our recent paper \cite{zhou17arxiv} we have carried out a comprehensive
study on the analysis and design of the complex-valued linear system
(\ref{sys}), including state response, controllability, observability,
stability, stabilization, pole assignment, linear quadratic regulation, and
state observer design. The obtained results will reduce to classical ones when
they are applied on the normal linear system (\ref{normal}), and will reduce
to and/or improve the existing results when they are applied on the antilinear
system (\ref{antilinear}). Our study is mainly based on the properties of the
bimatrix $\{P_{1},P_{2}\} \in \{ \mathbf{C}^{n\times m},\mathbf{C}^{n\times
m}\}$, especially, its real representation
\begin{equation}
\left \{  P_{1},P_{2}\right \}  _{\circ}\triangleq \left[
\begin{array}
[c]{cc}%
\mathrm{Re}\left(  P_{1}+P_{2}\right)  & -\mathrm{Im}\left(  P_{1}%
+P_{2}\right) \\
\mathrm{Im}\left(  P_{1}-P_{2}\right)  & \mathrm{Re}\left(  P_{1}%
-P_{2}\right)
\end{array}
\right]  \in \mathbf{R}^{2n\times2m}, \label{real}%
\end{equation}
and complex-lifting%
\begin{equation}
\left \{  P_{1},P_{2}\right \}  _{\diamond}\triangleq \left[
\begin{array}
[c]{cc}%
P_{1} & P_{2}^{\#}\\
P_{2} & P_{1}^{\#}%
\end{array}
\right]  \in \mathbf{C}^{2n\times2m}. \label{lifting1}%
\end{equation}
Particularly, we shown that, for stabilization and pole assignment of system
(\ref{sys}), the so-called full state feedback%
\begin{equation}
u=\left \{  K_{1},K_{2}\right \}  x, \label{eqfeedback}%
\end{equation}
is necessary \cite{zhou17arxiv}, and, generally, the well-known normal linear
feedback%
\begin{equation}
u=K_{1}x, \label{normalfeedback}%
\end{equation}
is valid only for the antilinear system (\ref{anticoprime}) when
$t\in \mathbf{Z}^{+}.$

In this paper, we continue to study the complex-valued linear system
(\ref{sys}). We are interested in the particular problem of pole assignment of
this class of systems by the full state feedback (\ref{eqfeedback}). We will
show that solutions to the pole assignment problem can be completely
characterized by solutions to a class of generalized Sylvester bimatrix
equations. Our study is clearly motivated by the existing work on pole
assignment of the normal linear system, for which it is well known that
solutions to the associated pole assignment can be characterized by solutions
to some (generalized) Sylvester matrix equations
\cite{bs82scl,duan93tac,duan15book}. We will provide complete solutions to
such a type of generalized Sylvester bimatrix equations and will show that the
obtained results include the existing ones for both the normal linear system
(\ref{normal}) and the antilinear system (\ref{antilinear}) as special cases.
From this point of view, we have built a quite general framework for pole
assignment of linear systems.

\subsection{\label{sec3.3}Second-Order Linear Systems}

In this subsection, we use a second-order linear system model to demonstrate
the purpose of studying the complex-valued linear system (\ref{sys}). Consider
the second-order linear system%
\begin{equation}
M\ddot{\xi}+D\dot{\xi}+K\xi=Gv, \label{second}%
\end{equation}
where $M,D,K\in \mathbf{R}^{n\times n}$ and $G\in \mathbf{R}^{n\times q}$ are
given matrices, $\xi$ is the state (often denotes the displacements of the
object to be controlled), and $v$ is the control. Let the initial condition be
$\xi \left(  0\right)  =\xi_{10}\in \mathbf{R}^{n}$ and $\dot{\xi}\left(
0\right)  =\xi_{20}\in \mathbf{R}^{n}.$ For simplicity, we only consider the
continuous-time case without output equation and assume that $M$ is
nonsingular. We further assume, without loss of generality, that $q=2m,$
since, otherwise, we set $G=[G,0_{n\times1}]$. Denote $G=[G_{1},G_{2}%
],G_{i}\in \mathbf{R}^{n\times m}$ and $v=[v_{1}^{\mathrm{T}},v_{2}%
^{\mathrm{T}}]^{\mathrm{T}},v_{i}\in \mathbf{R}^{m},i=1,2.$ Second-order linear
system can be used to describe many physical systems, for example, the
mass-spring system \cite{kailath80book}, and the spacecraft rendezvous control
system \cite{cw60jas}.

To describe the second-order linear system (\ref{second}) as a complex-valued
linear system, we first write it equivalently as
\begin{equation}
\left[
\begin{array}
[c]{c}%
\dot{\xi}\\
\ddot{\xi}%
\end{array}
\right]  =\left[
\begin{array}
[c]{cc}%
0_{n\times n} & I_{n}\\
-M^{-1}K & -M^{-1}D
\end{array}
\right]  \left[
\begin{array}
[c]{c}%
\xi \\
\dot{\xi}%
\end{array}
\right]  +\left[
\begin{array}
[c]{cc}%
0_{n\times m} & 0_{n\times m}\\
G_{1} & G_{2}%
\end{array}
\right]  \left[
\begin{array}
[c]{c}%
v_{1}\\
v_{2}%
\end{array}
\right]  . \label{second1}%
\end{equation}
If we choose%
\begin{equation}
x=\xi+\mathrm{j}\dot{\xi},\;u=v_{1}+\mathrm{j}v_{2}, \label{eqxu}%
\end{equation}
then system (\ref{second1}) can be equivalently
rewritten as (\ref{sys}) where the initial condition is $x_{0}=\xi
_{10}+\mathrm{j}\xi_{20},$ and $A_{i}\in \mathbf{C}^{n\times n},B_{i}%
\in \mathbf{C}^{n\times m},i=1,2,$ are given by
\begin{equation}
\left \{
\begin{array}
[c]{l}%
A_{1}=-\frac{1}{2}M^{-1}D-\frac{\mathrm{j}}{2}\left(  I_{n}+M^{-1}K\right)
,\\
A_{2}=\frac{1}{2}M^{-1}D-\frac{\mathrm{j}}{2}\left(  I_{n}-M^{-1}K\right)  ,\\
B_{1}=\frac{1}{2}G_{2}+\frac{\mathrm{j}}{2}G_{1},\\
B_{2}=-\frac{1}{2}G_{2}-\frac{\mathrm{j}}{2}G_{1}=-B_{1}.
\end{array}
\right.  \label{a12b12}%
\end{equation}

\begin{remark}
We have another method to describe system (\ref{second}) as (\ref{sys}). Write
(\ref{second}) as%
\begin{equation}
\left[
\begin{array}
[c]{c}%
\dot{\xi}\\
\ddot{\xi}%
\end{array}
\right]  =\left[
\begin{array}
[c]{cc}%
0_{n\times n} & I_{n}\\
-M^{-1}K & -M^{-1}D
\end{array}
\right]  \left[
\begin{array}
[c]{c}%
\xi \\
\dot{\xi}%
\end{array}
\right]  +\left[
\begin{array}
[c]{cc}%
0_{n\times q} & 0_{n\times q}\\
G & 0_{n\times q}%
\end{array}
\right]  \left[
\begin{array}
[c]{c}%
v\\
w
\end{array}
\right]  , \label{eq92}%
\end{equation}
where $w\in \mathbf{C}^{q}$ is a temp variable. Then, similar to (\ref{second1}%
), by setting $x$ as in (\ref{eqxu}) and $u=v+\mathrm{j}w,$ (\ref{eq92}) can
be written as (\ref{sys}), where $A_{i},i=1,2,$ are given by (\ref{a12b12})
and%
\[
B_{1}=\frac{\mathrm{j}}{2}G,\quad B_{2}=-\frac{\mathrm{j}}{2}G=-B_{1}.
\]
This method does not require that $q$ is an even number, but however leads to
higher dimensions of the inputs than (\ref{a12b12}).
\end{remark}

The corresponding complex-valued linear system model (\ref{sys}) for
(\ref{second}) seems more convenient to use than the augmented normal linear
system model (\ref{second1}) since it possesses the same dimension as the
original system (\ref{second}). We also mention that, as has been made clear
in \cite{zhou17arxiv}, the full state feedback (\ref{eqfeedback}) for the
associated complex-valued linear system (\ref{sys}) can be equivalently
written as
\begin{equation}
v=\left[
\begin{array}
[c]{c}%
v_{1}\\
v_{2}%
\end{array}
\right]  =\left \{  K_{1},K_{2}\right \}  _{\circ}\left[
\begin{array}
[c]{c}%
\xi \\
\dot{\xi}%
\end{array}
\right]  , \label{eqv}%
\end{equation}
where $\left \{  K_{1},K_{2}\right \}  _{\circ}$ is a real matrix. Therefore,
the full state feedback (\ref{eqfeedback}) can be implemented physically.

\subsection{Derivation of Linear Bimatrix Equations}

The problem of pole assignment for the complex-valued linear system
(\ref{sys}) by the full state feedback (\ref{eqfeedback}) can be stated as
finding the bimatrix $\left \{  K_{1},K_{2}\right \}  $\ $\in \{ \mathbf{C}%
^{m\times n},\mathbf{C}^{m\times n}\}$ such that the resulting closed-loop
system%
\begin{equation}
x^{+}=\left(  \left \{  A_{1},A_{2}\right \}  +\left \{  B_{1},B_{2}\right \}
\left \{  K_{1},K_{2}\right \}  \right)  x, \label{closed2}%
\end{equation}
possesses a desired eigenvalue set $\mathit{\Gamma}$ that is symmetric with
respect to the real axis. Since $\mathit{\Gamma}$ is symmetric with respect to
the real axis, there is a real matrix $F\in \mathbf{R}^{2n\times2n}$ such that
$\lambda \left(  F\right)  =\mathit{\Gamma}$.

\begin{lemma}
The pole assignment problem is solvable if and only if there exists a
nonsingular bimatrix $\{X_{1},X_{2}\} \in \{ \mathbf{C}^{n\times n}%
,\mathbf{C}^{n\times n}\}$ to the following generalized Sylvester bimatrix
equation%
\begin{equation}
\left \{  A_{1},A_{2}\right \}  \left \{  X_{1},X_{2}\right \}  +\left \{
B_{1},B_{2}\right \}  \left \{  Y_{1},Y_{2}\right \}  =\left \{  X_{1}%
,X_{2}\right \}  \left \{  F_{1},F_{2}\right \}  , \label{syl}%
\end{equation}
where $\left \{  F_{1},F_{2}\right \}  \in \{ \mathbf{C}^{n\times n}%
,\mathbf{C}^{n\times n}\}$ \ is the unique bimatrix satisfying
\begin{equation}
F=\left \{  F_{1},F_{2}\right \}  _{\circ}. \label{eqff}%
\end{equation}
In this case, the feedback gain bimatrix $\left \{  K_{1},K_{2}\right \}  $ is
determined by%
\begin{equation}
\left \{  K_{1},K_{2}\right \}  =\left \{  Y_{1},Y_{2}\right \}  \left \{
X_{1},X_{2}\right \}  ^{-1}. \label{eqsylk}%
\end{equation}

\end{lemma}

We give a remark regarding the closed-loop system (\ref{closed2}).

\begin{remark}
By the state transformation $y=\left \{  X_{1},X_{2}\right \}  ^{-1}x,$ we have
from (\ref{syl}) and (\ref{eqsylk}) that the closed-loop system (\ref{closed2}%
) is equivalent to%
\begin{align}
y^{+}  &  =\left \{  X_{1},X_{2}\right \}  ^{-1}x^{+}\nonumber \\
&  =\left \{  X_{1},X_{2}\right \}  ^{-1}\left(  \left \{  A_{1},A_{2}\right \}
+\left \{  B_{1},B_{2}\right \}  \left \{  K_{1},K_{2}\right \}  \right)
x\nonumber \\
&  =\left \{  X_{1},X_{2}\right \}  ^{-1}\left(  \left \{  A_{1},A_{2}\right \}
+\left \{  B_{1},B_{2}\right \}  \left \{  K_{1},K_{2}\right \}  \right)  \left \{
X_{1},X_{2}\right \}  y\nonumber \\
&  =\left \{  X_{1},X_{2}\right \}  ^{-1}\left(  \left \{  A_{1},A_{2}\right \}
\left \{  X_{1},X_{2}\right \}  +\left \{  B_{1},B_{2}\right \}  \left \{
Y_{1},Y_{2}\right \}  \right)  y\nonumber \\
&  =\left \{  X_{1},X_{2}\right \}  ^{-1}\left \{  X_{1},X_{2}\right \}  \left \{
F_{1},F_{2}\right \}  y\nonumber \\
&  =\left \{  F_{1},F_{2}\right \}  y. \label{closed1}%
\end{align}
Clearly, if we want the closed-loop system to be equivalent to a normal linear
system, then we should set $F_{2}=0_{n\times n},$ and to be equivalent to an
antilinear system, then we should set $F_{1}=0_{n\times n}.$
\end{remark}

Therefore, to solve the pole assignment problem for system (\ref{sys}), the
main task is to finding complete solutions to the generalized Sylvester
bimatrix equation (\ref{syl}), which is one of the main tasks in this paper
and will be studied in Sections \ref{sec3}-\ref{sec4}.

The generalized Sylvester bimatrix equation (\ref{syl}) is homogeneous and
thus its solution is non-unique. As done in pole assignment for normal linear
system \cite{bs82scl}, sometimes we may first prescribe $\{Y_{1},Y_{2}\}$ and
then seek the (unique) solution for $\{X_{1},X_{2}\}$ (if exists). In this
case, (\ref{syl}) becomes the following non-homogeneous one:%
\begin{equation}
\left \{  A_{1},A_{2}\right \}  \left \{  X_{1},X_{2}\right \}  -\left \{
X_{1},X_{2}\right \}  \left \{  F_{1},F_{2}\right \}  =\left \{  C_{1}%
,C_{2}\right \}  , \label{bimatrixsyl}%
\end{equation}
where $\{C_{1},C_{2}\} \triangleq-\left \{  B_{1},B_{2}\right \}  \left \{
Y_{1},Y_{2}\right \}  $ is known.\ This equation is referred to as the
Sylvester bimatrix equation, and will be studied in Section \ref{sec5}. The
non-homogeneous equation (\ref{bimatrixsyl}) also appears in the stability
analysis of the complex-valued linear system (\ref{sys}) with $t\in
\mathbf{R}^{+}$, which is asymptotically stable if and only if
\cite{zhou17arxiv}%
\begin{equation}
\left \{  A_{1},A_{2}\right \}  ^{\mathrm{H}}\left \{  P_{1},P_{2}\right \}
+\left \{  P_{1},P_{2}\right \}  \left \{  A_{1},A_{2}\right \}  =-\left \{
Q_{1},Q_{2}\right \}  , \label{bilya}%
\end{equation}
has a (unique) solution $\left \{  P_{1},P_{2}\right \}  >0$ for any given
$\left \{  Q_{1},Q_{2}\right \}  >0.$ Clearly, (\ref{bilya}) is in the form of
(\ref{bimatrixsyl}), and will be referred to as Lyapunov bimatrix equation.

When we study stability of the complex-valued linear system (\ref{sys}) with
$t\in \mathbf{Z}^{+},$ the discrete-time Lyapunov bimatrix equation%
\begin{equation}
\left \{  P_{1},P_{2}\right \}  =\left \{  A_{1},A_{2}\right \}  ^{\mathrm{H}%
}\left \{  P_{1},P_{2}\right \}  \left \{  A_{1},A_{2}\right \}  +\left \{
Q_{1},Q_{2}\right \}  , \label{bidilya}%
\end{equation}
is encountered. It is shown in \cite{zhou17arxiv} that stability of
(\ref{sys}) with $t\in \mathbf{Z}^{+}$ is equivalent to the existence of a
(unique) solution $\left \{  P_{1},P_{2}\right \}  >0$ to (\ref{bidilya}) for
any given $\left \{  Q_{1},Q_{2}\right \}  >0.$ Equation (\ref{bidilya}) is a
special case of%
\begin{equation}
\left \{  X_{1},X_{2}\right \}  =\left \{  A_{1},A_{2}\right \}  \left \{
X_{1},X_{2}\right \}  \left \{  F_{1},F_{2}\right \}  +\left \{  C_{1}%
,C_{2}\right \}  , \label{bistein}%
\end{equation}
which is referred to as the Stein bimatrix equation, and will also be studied
in Section \ref{sec5}.

Though in the above the bimatrix $\left \{  F_{1},F_{2}\right \}  $ has the same
dimension as $\left \{  A_{1},A_{2}\right \}  ,$ however, this is not necessary.
Hence, without loss of generality, hereafter we assume, if not specified, that
$\left \{  F_{1},F_{2}\right \}  \in \{ \mathbf{C}^{p\times p},\mathbf{C}%
^{p\times p}\},$ where $p$ is any positive integer.

\section{\label{sec3}Solutions to Generalized Sylvester Bimatrix Equations}

In this section we study solutions to the generalized Sylvester bimatrix
equation (\ref{syl}) and will also consider a special case that its
coefficients are determined by the normal linear system (\ref{normal}).
Hereafter we assume that system (\ref{sys}) (or $(\{A_{1},A_{2}\},\{B_{1}%
,B_{2}\})$) is controllable (see \cite{zhou17arxiv} for definition and
criterion for the controllability of system (\ref{sys})).

\subsection{General Solutions}

Two polynomial bimatrices $\left \{  N_{1}(s),N_{2}(s)\right \}  \in \{
\mathbf{C}^{n\times m},\mathbf{C}^{n\times m}\}$ and $\left \{  D_{1}%
(s),D_{2}(s)\right \}  \in \{ \mathbf{C}^{m\times m},\mathbf{C}^{m\times m}\}$
are said to be right-coprime if 
\begin{equation}
2m=\mathrm{rank}\left \{  \left[
\begin{array}
[c]{c}%
N_{1}\left(  s\right) \\
D_{1}\left(  s\right)
\end{array}
\right]  ,\left[
\begin{array}
[c]{c}%
N_{2}\left(  s\right) \\
D_{2}\left(  s\right)
\end{array}
\right]  \right \}  ,\forall s\in \mathbf{C,} \label{coprime1}%
\end{equation}
where, and hereafter, $s$ should be treated as a \textit{real parameter} when
computing its conjugate, namely, $s=s^{\#}$. We then can present the following
explicit solutions to the generalized Sylvester bimatrix equation (\ref{syl}).

\begin{theorem}
\label{th1}Let $\left \{  N_{1}(s),N_{2}(s)\right \}  \in \{ \mathbf{C}^{n\times
m},\mathbf{C}^{n\times m}\}$ and $\left \{  D_{1}(s),D_{2}(s)\right \}  \in \{
\mathbf{C}^{m\times m},\mathbf{C}^{m\times m}\}$ be two polynomial bimatrices
such that%
\begin{equation}
\left \{  sI_{n}-A_{1},-A_{2}\right \}  \left \{  N_{1}(s),N_{2}(s)\right \}
=\left \{  B_{1},B_{2}\right \}  \left \{  D_{1}(s),D_{2}(s)\right \}  ,
\label{coprime}%
\end{equation}
and $\left \{  N_{1}(s),N_{2}(s)\right \}  $ and $\left \{  D_{1}(s),D_{2}%
(s)\right \}  $ are right-coprime. Then complete solutions to (\ref{syl}) are
given by%
\begin{equation}
\left \{
\begin{array}
[c]{rl}%
\left \{  X_{1},X_{2}\right \}  & =\sum \limits_{i=0}^{\omega}\left \{
N_{1i},N_{2i}\right \}  \left \{  Z_{1},Z_{2}\right \}  \left \{  F_{1}%
,F_{2}\right \}  ^{i},\\
\left \{  Y_{1},Y_{2}\right \}  & =\sum \limits_{i=0}^{\omega}\left \{
D_{1i},D_{2i}\right \}  \left \{  Z_{1},Z_{2}\right \}  \left \{  F_{1}%
,F_{2}\right \}  ^{i},
\end{array}
\right.  \label{xysolution}%
\end{equation}
where $\{Z_{1},Z_{2}\} \in \{ \mathbf{C}^{m\times p},\mathbf{C}^{m\times p}\}$
is an arbitrarily bimatrix and%
\begin{equation}
\left[
\begin{array}
[c]{c}%
N_{j}(s)\\
D_{j}(s)
\end{array}
\right]  =\sum \limits_{i=0}^{\omega}\left[
\begin{array}
[c]{c}%
N_{ji}\\
D_{ji}%
\end{array}
\right]  s^{i},\;j=1,2,\; \omega \in \mathbf{Z}^{+}. \label{nidi}%
\end{equation}

\end{theorem}

Notice that the polynomial bimatrix equations (\ref{coprime}%
) can be written as coupled polynomial matrix equations%
\begin{equation}
\left \{
\begin{array}
[c]{rl}%
\left(  sI_{n}-A_{1}\right)  N_{1}(s)-A_{2}^{\#}N_{2}(s) & =B_{1}%
D_{1}(s)+B_{2}^{\#}D_{2}(s),\\
\left(  sI_{n}-A_{1}\right)  N_{2}^{\#}(s)-A_{2}^{\#}N_{1}^{\#}(s) &
=B_{1}D_{2}^{\#}(s)+B_{2}^{\#}D_{1}^{\#}(s),
\end{array}
\right.  \label{eqcoprime_coupled}%
\end{equation}
and, the generalized Sylvester bimatrix equation (\ref{syl}) can also be
expressed equivalently as coupled matrix equations%
\begin{equation}
\left \{
\begin{array}
[c]{rl}%
A_{1}X_{1}+A_{2}^{\#}X_{2}+B_{1}Y_{1}+B_{2}^{\#}Y_{2} & =X_{1}F_{1}+X_{2}%
^{\#}F_{2},\\
A_{1}X_{2}^{\#}+A_{2}^{\#}X_{1}^{\#}+B_{1}Y_{2}^{\#}+B_{2}^{\#}Y_{1}^{\#} &
=X_{1}F_{2}^{\#}+X_{2}^{\#}F_{1}^{\#}.
\end{array}
\right.  \label{sylcoupled}%
\end{equation}
In the next subsection, we show how to transform these coupled equations into
equivalent decoupled ones.

\subsection{Decoupling of Coupled Equations}

We first show that the coupled polynomial matrix equations in
(\ref{eqcoprime_coupled}) can be decoupled.

\begin{lemma}
\label{lm1}The coupled polynomial matrix equations (\ref{eqcoprime_coupled})
with unknowns $(N_{1}(s),N_{2}(s),D_{1}(s),D_{2}(s))$ are solvable if and only
if the following decoupled matrix equations%
\begin{equation}
\left \{
\begin{array}
[c]{rl}%
\left(  sI_{n}-A_{1}\right)  N_{+}(s)-A_{2}^{\#}N_{+}^{\#}(s) & =B_{1}%
D_{+}(s)+B_{2}^{\#}D_{+}^{\#}(s),\\
\left(  sI_{n}-A_{1}\right)  N_{-}(s)+A_{2}^{\#}N_{-}^{\#}(s) & =B_{1}%
D_{-}(s)-B_{2}^{\#}D_{-}^{\#}(s),
\end{array}
\right.  \label{eqcoprime_decoupled}%
\end{equation}
with unknowns $(N_{+}(s),N_{-}(s),D_{+}(s),D_{-}(s))$ are solvable. Moreover,
$(N_{1}(s),N_{2}(s),D_{1}(s),D_{2}(s))$ and $(N_{+}(s)$, $N_{-}(s),D_{+}%
(s),D_{-}(s))$ are one-to-one according to%
\begin{equation}
\left \{
\begin{array}
[c]{ll}%
N_{1}(s)=\frac{1}{2}\left(  N_{+}(s)+N_{-}(s)\right)  , & D_{1}(s)=\frac{1}%
{2}\left(  D_{+}(s)+D_{-}(s)\right)  ,\\
N_{2}(s)=\frac{1}{2}\left(  N_{+}(s)-N_{-}(s)\right)  ^{\#}, & D_{2}%
(s)=\frac{1}{2}\left(  D_{+}(s)-D_{-}(s)\right)  ^{\#}.
\end{array}
\right.  \label{eqnpnm}%
\end{equation}
Furthermore, $\left \{  N_{1}(s),N_{2}(s)\right \}  $ and $\{D_{1}%
(s),D_{2}(s)\}$ are right-coprime if and only if%
\begin{equation}
\mathrm{rank}\left[
\begin{array}
[c]{cc}%
N_{+}(s) & -N_{-}(s)\\
D_{+}(s) & -D_{-}(s)\\
N_{+}^{\#}(s) & N_{-}^{\#}(s)\\
D_{+}^{\#}(s) & D_{-}^{\#}(s)
\end{array}
\right]  =2m,\; \forall s\in \mathbf{C}. \label{eqcoprime5}%
\end{equation}

\end{lemma}

By this lemma, the two matrix pairs $\left(  N_{+}(s),D_{+}(s)\right)  $ and
$\left(  N_{-}(s),D_{-}(s)\right)  $ can be solved separately, which is useful
in computation. In fact, we need only to solve the first equation since the
second one can be solved in a similar way by setting $A_{2}\mapsto-A_{2}$ and
$B_{2}\mapsto-B_{2}.$ General solutions to (\ref{eqcoprime_decoupled}) and
(\ref{eqcoprime5}) will be investigated in a separate paper.

Similar to Lemma \ref{lm1}, the coupled matrix equations (\ref{sylcoupled})
can also be decoupled in some cases, as shown below.

\begin{lemma}
\label{lm2}Assume that there exist two real matrices $F_{ii}\in \mathbf{R}%
^{p\times p},i=1,2,$ such that
\begin{equation}
F\triangleq \{F_{1},F_{2}\}_{\circ}=\mathrm{diag}\{F_{11},F_{22}\}, \label{eqf}%
\end{equation}
namely, $F_{1}$ and $F_{2}$ are also real matrices given by (see (\ref{real}))%
\begin{equation}
F_{1}=\frac{1}{2}\left(  F_{11}+F_{22}\right)  ,\;F_{2}=\frac{1}{2}\left(
F_{11}-F_{22}\right)  . \label{eqf1f2}%
\end{equation}
Then the associated coupled matrix equations (\ref{sylcoupled}) are solvable
with unknowns $(X_{1},X_{2},Y_{1},Y_{2})$ if and only if the decoupled matrix
equations%
\begin{equation}
\left \{
\begin{array}
[c]{rl}%
A_{1}X_{+}+A_{2}^{\#}X_{+}^{\#}+B_{1}Y_{+}+B_{2}^{\#}Y_{+}^{\#} &
=X_{+}\left(  F_{1}+F_{2}\right)  ,\\
A_{1}X_{-}-A_{2}^{\#}X_{-}^{\#}+B_{1}Y_{-}-B_{2}^{\#}Y_{-}^{\#} &
=X_{-}\left(  F_{1}-F_{2}\right)  ,
\end{array}
\right.  \label{eqsylde}%
\end{equation}
are solvable with unknowns $(X_{+},Y_{+})$ and $(X_{-},Y_{-}).$ Moreover,
$(X_{1},X_{2},Y_{1},Y_{2})\ $and $(X_{+},Y_{+},X_{-},Y_{-})$ are one-to-one
according to%
\begin{equation}
\left \{
\begin{array}
[c]{ll}%
X_{1}=\frac{1}{2}\left(  X_{+}+X_{-}\right)  , & Y_{1}=\frac{1}{2}\left(
Y_{+}+Y_{-}\right)  ,\\
X_{2}=\frac{1}{2}\left(  X_{+}-X_{-}\right)  ^{\#}, & Y_{2}=\frac{1}{2}\left(
Y_{+}-Y_{-}\right)  ^{\#}.
\end{array}
\right.  \label{eqxpxm}%
\end{equation}

\end{lemma}

The assumption (\ref{eqf}) is not restrictive if only asymptotic stability of
the closed-loop system is concerned since, in view of (\ref{eqff}), we can
always find real matrices $F_{1}\ $and $F_{2}$ such that $F$ is asymptotically
stable. Combining Lemmas \ref{lm1} and \ref{lm2} gives the following result.

\begin{theorem}
Assume that there exist two real matrices $F_{ii}\in \mathbf{R}^{p\times
p},i=1,2$ satisfying (\ref{eqf}) and $\left(  F_{1},F_{2}\right)  $ is
given\ by (\ref{eqf1f2}). Let $(N_{+}(s),D_{+}(s))$ and $(N_{-}(s),D_{-}(s))$
satisfy respectively the first and second equations of
(\ref{eqcoprime_decoupled}) and such that (\ref{eqcoprime5}). Denote%
\begin{equation}
\left[
\begin{array}
[c]{c}%
N_{\pm}(s)\\
D_{\pm}(s)
\end{array}
\right]  =\sum \limits_{i=0}^{\omega}\left[
\begin{array}
[c]{c}%
N_{\pm i}\\
D_{\pm i}%
\end{array}
\right]  s^{i},\; \omega \in \mathbf{Z}^{+}. \label{ndpmpm}%
\end{equation}
Then complete solutions to (\ref{eqsylde}) are given by%
\begin{equation}
\left[
\begin{array}
[c]{c}%
X_{\pm}\\
Y_{\pm}%
\end{array}
\right]  =\frac{1}{2}\sum \limits_{i=0}^{\omega}\left(  \left[
\begin{array}
[c]{c}%
N_{\pm i}\\
D_{\pm i}%
\end{array}
\right]  \left(  Z_{\pm}+Z_{\pm}^{\#}\right)  +\left[
\begin{array}
[c]{c}%
N_{\mp i}\\
D_{\mp i}%
\end{array}
\right]  \left(  Z_{\pm}-Z_{\pm}^{\#}\right)  \right)  \left(  F_{1}\pm
F_{2}\right)  ^{i}, \label{eqxypm}%
\end{equation}
where $\left(  Z_{+},Z_{-}\right)  $ and $\left(  Z_{1},Z_{2}\right)  $ are
one-to-one according to%
\begin{equation}
Z_{\pm}=Z_{1}\pm Z_{2}^{\#}. \label{eqzpm}%
\end{equation}

\end{theorem}

It follows that, though $\left(  X_{+},Y_{+}\right)  $ and $\left(
X_{-},Y_{-}\right)  $ are decoupled in (\ref{eqsylde}), and $(N_{+}%
(s),D_{+}(s))$ and $(N_{-}(s),D_{-}(s))$ are also decoupled in
(\ref{eqcoprime_decoupled}), $\left(  X_{+},Y_{+}\right)  $ (and $\left(
X_{-},Y_{-}\right)  $) depends on both $(N_{+}(s),D_{+}(s))$ and
$(N_{-}(s),D_{-}(s)).$ If we let $Z_{\pm}\in \mathbf{R}^{m\times p},$ then they
are decoupled as%
\[
\left[
\begin{array}
[c]{c}%
X_{\pm}\\
Y_{\pm}%
\end{array}
\right]  =\sum \limits_{i=0}^{\omega}\left[
\begin{array}
[c]{c}%
N_{\pm i}\\
D_{\pm i}%
\end{array}
\right]  Z_{\pm}\left(  F_{1}\pm F_{2}\right)  ^{i},
\]
which is appealing in mathematics.

\subsection{\label{sec4.3}Solutions for Normal Linear Systems}

When considering the generalized Sylvester bimatrix equation (\ref{syl}) for
the normal linear system (\ref{normal}), we may want the closed-loop system to
be (similar to) a normal system as well (see Remark \ref{rm4} later for a
different situation). Thus, according to (\ref{closed1}), we should choose
$F_{2}=0_{n\times n}.$ Consequently, the generalized Sylvester bimatrix
equation (\ref{syl}) or the coupled matrix equation (\ref{sylcoupled})
becomes
\begin{equation}
\left \{
\begin{array}
[c]{rl}%
A_{1}X_{1}+B_{1}Y_{1} & =X_{1}F_{1},\\
A_{1}X_{2}^{\#}+B_{1}Y_{2}^{\#} & =X_{2}^{\#}F_{1}^{\#},
\end{array}
\right.  \label{sylnormal}%
\end{equation}
where $F_{1}$ is such that $y^{+}=F_{1}y$ is asymptotically stable. These two
equations in (\ref{sylnormal}) are exactly the same one%
\begin{equation}
A_{1}X_{0}+B_{1}Y_{0}=X_{0}F_{0}, \label{sylvester}%
\end{equation}
where $F_{0}=F_{1}$ and $F_{1}^{\#}.$ Equation (\ref{sylvester}) is known as
the generalized Sylvester matrix equation and has been extensively used and
studied in the literature for pole assignment of the normal linear system
(\ref{normal}) \cite{duan93tac,duan15book,zd06scl}.

In this case, the two polynomial matrix equations (\ref{eqcoprime_decoupled})
reduce further to the single one%
\begin{equation}
\left(  sI_{n}-A_{1}\right)  N_{0}(s)=B_{1}D_{0}\left(  s\right)  ,
\label{normalcoprime}%
\end{equation}
that has been investigated in the literature, for example, \cite{duan93tac}
and \cite{zdl09iet}. We can simply choose
\begin{equation}
\left \{
\begin{array}
[c]{ll}%
N_{1}(s)=N_{0}(s), & N_{2}(s)=0_{n\times m},\\
D_{1}(s)=D_{0}(s),\; & D_{2}(s)=0_{m\times m},
\end{array}
\right.  \label{eq82}%
\end{equation}
to satisfy (\ref{eqcoprime_coupled}). Thus (\ref{coprime1}) is equivalent to,
for all $s\in \mathbf{C},$%
\[
2m=\mathrm{rank}\left[
\begin{array}
[c]{cc}%
N_{1}\left(  s\right)  & 0_{n\times m}\\
D_{1}\left(  s\right)  & 0_{m\times m}\\
0_{n\times m} & N_{1}^{\#}\left(  s\right) \\
0_{m\times m} & D_{1}^{\#}\left(  s\right)
\end{array}
\right]  =2\mathrm{rank}\left[
\begin{array}
[c]{c}%
N_{0}(s)\\
D_{0}(s)
\end{array}
\right]  ,
\]
which implies that $\left \{  N_{1}(s),N_{2}(s)\right \}  $ and $\{D_{1}%
(s),D_{2}(s)\}$ are right-coprime if and only if
\[
\mathrm{rank}\left[
\begin{array}
[c]{c}%
N_{0}(s)\\
D_{0}(s)
\end{array}
\right]  =m,\; \forall s\in \mathbf{C},
\]
namely, $N_{0}\left(  s\right)  $ and $D_{0}\left(  s\right)  $ are
right-coprime in the normal sense \cite{kailath80book}. In this case, by
denoting%
\[
\left[
\begin{array}
[c]{c}%
N_{0}(s)\\
D_{0}(s)
\end{array}
\right]  =\sum \limits_{i=0}^{\omega}\left[
\begin{array}
[c]{c}%
N_{0i}\\
D_{0i}%
\end{array}
\right]  s^{i},
\]
complete solutions to (\ref{sylvester}) are given by \cite{zd06scl}%
\begin{equation}
\left[
\begin{array}
[c]{c}%
X_{0}\\
Y_{0}%
\end{array}
\right]  =\sum \limits_{i=0}^{\omega}\left[
\begin{array}
[c]{c}%
N_{0i}\\
D_{0i}%
\end{array}
\right]  Z_{0}F_{0}^{i}, \label{eq91}%
\end{equation}
where $Z_{0}\in \mathbf{C}^{m\times p}$ is any matrix. On the other hand, in
view of (\ref{eq82}), we can apply solution (\ref{xysolution}) on equation
(\ref{sylnormal}) (setting $A_{2}=0_{n\times n},B_{2}=0_{n\times m}%
,F_{2}=0_{p\times p}$ and using (\ref{eq82})) to obtain%
\begin{equation}
\left \{
\begin{array}
[c]{l}%
\left[
\begin{array}
[c]{c}%
X_{1}\\
Y_{1}%
\end{array}
\right]  =\sum \limits_{i=0}^{\omega}\left[
\begin{array}
[c]{c}%
N_{0i}\\
D_{0i}%
\end{array}
\right]  Z_{1}F_{1}^{i},\\
\left[
\begin{array}
[c]{c}%
X_{2}\\
Y_{2}%
\end{array}
\right]  =\sum \limits_{i=0}^{\omega}\left[
\begin{array}
[c]{c}%
N_{0i}\\
D_{0i}%
\end{array}
\right]  Z_{2}F_{1}^{\#i},
\end{array}
\right.  \label{eqxyxy}%
\end{equation}
where $Z_{i}\in \mathbf{C}^{m\times p},i=1,2,$ are any matrices. These two
expressions coincide exactly with (\ref{eq91}).

\begin{remark}
\label{rm9}Notice that, if we choose $Z_{2}=0_{m\times n},$ then we have from
(\ref{eqxyxy}) that $X_{2}=0_{n\times n}$ and $Y_{2}=0_{m\times n}$. Then, if
$Z_{1}$ is properly chosen such that $X_{1}$ is nonsingular, by (\ref{eqsylk}%
), the feedback gain bimatrix is given by%
\begin{align}
\left \{  K_{1},K_{2}\right \}   &  =\left \{  Y_{1},0_{m\times n}\right \}
\left \{  X_{1},0_{n\times n}\right \}  ^{-1}\nonumber \\
&  =\left \{  Y_{1},0_{m\times n}\right \}  \left \{  X_{1}^{-1},0_{n\times
n}\right \} \nonumber \\
&  =\left \{  Y_{1}X_{1}^{-1},0_{m\times n}\right \}  , \label{eq89}%
\end{align}
which means that the resulting controller is just the normal linear state
feedback (\ref{normalfeedback}). However, if $Z_{2}\neq0_{m\times n},$ then
$X_{2}\neq0_{n\times n}$ and $Y_{2}\neq0_{m\times n},$ which in turn implies
that the resulting controller is the full state feedback (\ref{eqfeedback}).
Yet in both cases the closed-loop system is equivalent to a normal linear
system (see (\ref{closed1})).
\end{remark}

\section{\label{sec4}Solutions for Antilinear Systems}

In this section, we carry out a careful study on pole assignment for the
antilinear system (\ref{antilinear}) and present explicit solutions to the
associated Sylvester bimatrix equations (\ref{syl}) in different special cases.

\subsection{General Solutions}

We first provide a method for computing the right-coprime factorization
(\ref{coprime}) or (\ref{eqcoprime_decoupled}) associated with the antilinear
system (\ref{antilinear}).

\begin{lemma}
\label{coro10}Consider the antilinear system (\ref{antilinear}). Let $\left(
N_{0}(s),D_{0}\left(  s\right)  \right)  \in(\mathbf{C}^{n\times m}%
,\mathbf{C}^{m\times m})$ satisfy
\begin{equation}
sN_{0}^{\#}(s)-A_{2}N_{0}(s)=B_{2}D_{0}(s). \label{anticoprime}%
\end{equation}
Then the polynomial matrices $N_{+}(s),D_{+}(s),N_{-}(s)$ and $D_{-}(s)\}$
satisfying (\ref{eqcoprime_decoupled}) can be chosen as
\begin{equation}
\left \{
\begin{array}
[c]{ll}%
N_{+}\left(  s\right)  =N_{0}\left(  s\right)  , & D_{+}\left(  s\right)
=D_{0}\left(  s\right)  ,\\
N_{-}\left(  s\right)  =N_{0}\left(  -s\right)  , & D_{-}\left(  s\right)
=D_{0}\left(  -s\right)  ,
\end{array}
\right.  \label{eq76}%
\end{equation}
or equivalently, $\{N_{1}(s),N_{2}(s)\}$ and $\{D_{1}(s),D_{2}(s)\}$
satisfying (\ref{coprime}) can be chosen as%
\begin{equation}
\left \{
\begin{array}
[c]{ll}%
N_{1}\left(  s\right)  =\frac{1}{2}\left(  N_{0}\left(  s\right)
+N_{0}\left(  -s\right)  \right)  , & D_{1}(s)=\frac{1}{2}\left(
D_{0}(s)+D_{0}(-s)\right)  ,\\
N_{2}\left(  s\right)  =\frac{1}{2}\left(  N_{0}\left(  s\right)
-N_{0}\left(  -s\right)  \right)  ^{\#}, & D_{2}(s)=\frac{1}{2}\left(
D_{0}(s)-D_{0}(-s)\right)  ^{\#}.
\end{array}
\right.  \label{eqcoprime6}%
\end{equation}
Moreover, $\{N_{1}(s),N_{2}(s)\}$ and $\{D_{1}(s),D_{2}(s)\}$ are
right-coprime if and only if%
\begin{equation}
\mathrm{rank}\left[
\begin{array}
[c]{cc}%
N_{0}\left(  s\right)  & -N_{0}\left(  -s\right) \\
D_{0}\left(  s\right)  & -D_{0}\left(  -s\right) \\
N_{0}^{\#}\left(  s\right)  & N_{0}^{\#}\left(  -s\right) \\
D_{0}^{\#}\left(  s\right)  & D_{0}^{\#}\left(  -s\right)
\end{array}
\right]  =2m,\; \forall s\in \mathbf{C}. \label{eq79}%
\end{equation}

\end{lemma}

A polynomial matrix pair $\left(  N_{0}(s),D_{0}\left(  s\right)  \right)  $
satisfying (\ref{eq79}) may be called anti-right-coprime, and the equation
(\ref{anticoprime}) can be referred to as the anti-right-coprime factorization
of the antilinear system (\ref{antilinear}), which can be studied by using the
approach in \cite{zdl09iet}.

The Sylvester bimatrix equation (\ref{syl}) or the decoupled matrix equations
(\ref{eqsylde}) associated with the antilinear system (\ref{antilinear}) is
equivalent to%
\begin{equation}
\left \{
\begin{array}
[c]{cc}%
A_{2}^{\#}X_{+}^{\#}+B_{2}^{\#}Y_{+}^{\#} & =X_{+}\left(  F_{1}+F_{2}\right)
,\\
-A_{2}^{\#}X_{-}^{\#}-B_{2}^{\#}Y_{-}^{\#} & =X_{-}\left(  F_{1}-F_{2}\right)
.
\end{array}
\right.  \label{eqsyl2}%
\end{equation}
These two equations take also the same form while their coefficients are different.

\begin{corollary}
\label{coro14}Assume that there exist two real matrices $F_{ii}\in
\mathbf{R}^{p\times p},i=1,2$ satisfying (\ref{eqf}) and $\left(  F_{1}%
,F_{2}\right)  $ is given\ by (\ref{eqf1f2}). Let $\left(  N_{0}%
(s),D_{0}\left(  s\right)  \right)  \in(\mathbf{C}^{n\times m},\mathbf{C}%
^{m\times m})$ satisfy (\ref{anticoprime}) and (\ref{eq79}). Denote%
\begin{equation}
\left[
\begin{array}
[c]{c}%
N_{0}\left(  s\right) \\
D_{0}\left(  s\right)
\end{array}
\right]  =\sum \limits_{i=0}^{\omega}\left[
\begin{array}
[c]{c}%
N_{0,i}\\
D_{0,i}%
\end{array}
\right]  s^{i}\triangleq \sum \limits_{i=0}^{2\varpi}\left[
\begin{array}
[c]{c}%
N_{0,i}\\
D_{0,i}%
\end{array}
\right]  s^{i}. \label{eqn0d0}%
\end{equation}
Then complete solutions to (\ref{eqsyl2}) are given by%
\begin{equation}
\left[
\begin{array}
[c]{c}%
X_{\pm}\\
Y_{\pm}%
\end{array}
\right]  =\sum \limits_{i=0}^{\varpi}\left[
\begin{array}
[c]{c}%
N_{0,2i}\\
D_{0,2i}%
\end{array}
\right]  Z_{\pm}\left(  F_{1}\pm F_{2}\right)  ^{2i}\pm \sum \limits_{i=0}%
^{\varpi-1}\left[
\begin{array}
[c]{c}%
N_{0,2i+1}\\
D_{0,2i+1}%
\end{array}
\right]  Z_{\pm}^{\#}\left(  F_{1}\pm F_{2}\right)  ^{2i+1}, \label{eq99}%
\end{equation}
where $Z_{\pm}$ are determined by (\ref{eqzpm}).
\end{corollary}

\subsection{Normalization of Antilinear Systems}

We now consider a special case that $F_{2}=0_{n\times n}.$ In this case, by
(\ref{closed1}), the closed-loop system (\ref{closed2}) is equivalent to
$y^{+}=F_{1}y,$ which is a normal linear system and is asymptotically stable
if and only if $F_{1}$ is asymptotically stable ($\mu \left(  F_{1}\right)  <0$
when $t\in \mathbf{R}^{+}$ and $\rho \left(  F_{1}\right)  <1$ when
$t\in \mathbf{Z}^{+}).$ In this case, the generalized Sylvester bimatrix
equation (\ref{syl}) or the coupled matrix equations (\ref{sylcoupled})
become
\begin{equation}
\left \{
\begin{array}
[c]{rl}%
A_{2}^{\#}X_{2}+B_{2}^{\#}Y_{2} & =X_{1}F_{1},\\
A_{2}^{\#}X_{1}^{\#}+B_{2}^{\#}Y_{1}^{\#} & =X_{2}^{\#}F_{1}^{\#}.
\end{array}
\right.  \label{conjsyl2}%
\end{equation}
Hence, in this case, if (\ref{conjsyl2}) has a nonsingular solution
$\{X_{1},X_{2}\},$ then the full state feedback (\ref{eqfeedback}) and
(\ref{eqsylk}) will make the system be a normal linear system. We call such a
procedure as the \textit{normalization} of the antilinear system
(\ref{antilinear}).

\begin{corollary}
Let $F_{2}=0_{p\times p}$ and $\left(  N_{0}(s),D_{0}\left(  s\right)
\right)  \in(\mathbf{C}^{n\times m},\mathbf{C}^{m\times m})$ satisfy
(\ref{anticoprime}), (\ref{eq79}) and (\ref{eqn0d0}). Then complete solutions
to the first equation of (\ref{conjsyl2}) are given by%
\begin{equation}
\left[
\begin{array}
[c]{c}%
X_{1}\\
Y_{1}%
\end{array}
\right]  =\sum \limits_{i=0}^{\omega}\left \{
\begin{array}
[c]{ll}%
\left[
\begin{array}
[c]{c}%
N_{0,i}\\
D_{0,i}%
\end{array}
\right]  Z_{1}F_{1}^{i}, & i\text{ is even}\\
\left[
\begin{array}
[c]{c}%
N_{0,i}\\
D_{0,i}%
\end{array}
\right]  Z_{2}F_{1}^{i}, & i\text{ is odd,}%
\end{array}
\right.  \label{x1y1a}%
\end{equation}
and complete solutions to the second equation of (\ref{conjsyl2}) are given by%
\begin{equation}
\left[
\begin{array}
[c]{c}%
X_{2}\\
Y_{2}%
\end{array}
\right]  =\sum \limits_{i=0}^{\omega}\left \{
\begin{array}
[c]{ll}%
\left[
\begin{array}
[c]{c}%
N_{0,i}^{\#}\\
D_{0,i}^{\#}%
\end{array}
\right]  Z_{2}F_{1}^{i}, & i\text{ is even}\\
\left[
\begin{array}
[c]{c}%
N_{0,i}^{\#}\\
D_{0,i}^{\#}%
\end{array}
\right]  Z_{1}F_{1}^{i}, & i\text{ is odd,}%
\end{array}
\right.  \label{x2y2a}%
\end{equation}
\newline where $Z_{i}\in \mathbf{C}^{m\times p},i=1,2$ are arbitrary matrices.
\end{corollary}

\subsection{Anti-Preserving of Antilinear Systems}

We next consider a special case that $F_{1}=0_{n\times n}.$ In this case, by
(\ref{closed1}), the closed-loop system (\ref{closed2}) is equivalent to
$y^{+}=F_{2}^{\#}y^{\#},$ which is still an antilinear system and is
asymptotically stable if and only if $t\in \mathbf{Z}^{+}$ and
\cite{zhou17arxiv}
\begin{equation}
\rho \left(  F_{2}F_{2}^{\#}\right)  <1. \label{eqf2rho}%
\end{equation}
In this case, the generalized Sylvester bimatrix equation (\ref{syl}) or the
coupled matrix equations (\ref{sylcoupled}) become
\begin{equation}
\left \{
\begin{array}
[c]{rl}%
A_{2}^{\#}X_{1}^{\#}+B_{2}^{\#}Y_{1}^{\#} & =X_{1}F_{2}^{\#},\\
A_{2}^{\#}X_{2}+B_{2}^{\#}Y_{2} & =X_{2}^{\#}F_{2},
\end{array}
\right.  \label{conjsyl1}%
\end{equation}
which are decoupled. Hence, in this case, if (\ref{conjsyl1}) has a
nonsingular solution $\{X_{1},X_{2}\},$ then the full state feedback
(\ref{eqfeedback}) and (\ref{eqsylk}) will make the system be (equivalent to)
an antilinear system as well. We call such a procedure as the
\textit{anti-preserving} of the antilinear system (\ref{antilinear}).

\begin{corollary}
\label{coro15}Let $F_{1}=0_{p\times p}$ and $\left(  N_{0}(s),D_{0}\left(
s\right)  \right)  \in(\mathbf{C}^{n\times m},\mathbf{C}^{m\times m})$ satisfy
(\ref{anticoprime}), (\ref{eq79}) and (\ref{eqn0d0}). Then complete solutions
to the first equation of (\ref{conjsyl1}) are given by%
\begin{equation}
\left[
\begin{array}
[c]{c}%
X_{1}\\
Y_{1}%
\end{array}
\right]  =\sum \limits_{i=0}^{\omega}\left \{
\begin{array}
[c]{ll}%
\left[
\begin{array}
[c]{c}%
N_{0,i}\\
D_{0,i}%
\end{array}
\right]  Z_{1}\left(  F_{2}^{\#}F_{2}\right)  ^{\frac{i}{2}}, & i\text{ is
even}\\
\left[
\begin{array}
[c]{c}%
N_{0,i}\\
D_{0,i}%
\end{array}
\right]  Z_{1}^{\#}F_{2}\left(  F_{2}^{\#}F_{2}\right)  ^{\frac{i-1}{2}}, &
i\text{ is odd,}%
\end{array}
\right.  \label{eqx1y1}%
\end{equation}
and complete solutions to the second equation of (\ref{conjsyl1}) are given by%
\begin{equation}
\left[
\begin{array}
[c]{c}%
X_{2}\\
Y_{2}%
\end{array}
\right]  =\sum \limits_{i=0}^{\omega}\left \{
\begin{array}
[c]{ll}%
\left[
\begin{array}
[c]{c}%
N_{0,i}^{\#}\\
D_{0,i}^{\#}%
\end{array}
\right]  Z_{2}\left(  F_{2}^{\#}F_{2}\right)  ^{\frac{i}{2}}, & i\text{ is
even}\\
\left[
\begin{array}
[c]{c}%
N_{0,i}^{\#}\\
D_{0,i}^{\#}%
\end{array}
\right]  Z_{2}^{\#}F_{2}\left(  F_{2}^{\#}F_{2}\right)  ^{\frac{i-1}{2}}, &
i\text{ is odd,}%
\end{array}
\right.  \label{eqx2y2}%
\end{equation}
\newline where $Z_{i}\in \mathbf{C}^{m\times p},i=1,2$ are arbitrary matrices.
\end{corollary}

The solution $\left(  X_{1},Y_{1}\right)  $ given in Corollary \ref{coro15}
coincides with those obtained in \cite{wz17book}. We emphasize that the
solutions in Corollary \ref{coro15} can \textit{only be used to design
discrete-time antilinear systems}. In contrast, the coupled equations
(\ref{eqsyl2}) and (\ref{conjsyl2}) can be used to design both continuous-time
and discrete-time antilinear systems.

\begin{remark}
If we choose $Z_{2}=0_{m\times n}$ and $Z_{1}$ such that $X_{1}$ is
nonsingular, then, by (\ref{eqsylk}), we also have (\ref{eq89}), namely, the
resulting controller is the normal state feedback (\ref{normalfeedback}). This
case is just the one studied in \cite{wz17book} and our controller is exactly
the one obtained there (yet the rank condition (\ref{eq79}) was not available
in \cite{wz17book} where a different concept was adopted). However, similar to
Remark \ref{rm9}, if we choose $Z_{2}\neq0_{m\times n},$ such that $X_{2}%
\neq0_{n\times n}$ and/or $Y_{2}\neq0_{m\times n}$, then the resulting
controller is the full state feedback (\ref{eqfeedback}).
\end{remark}

\begin{remark}
\label{rm4}In Subsection \ref{sec4.3} we have assumed that $F_{2}=0_{n\times
n}$ for the normal linear system (\ref{normal}), namely, the closed-loop
system is equivalent to a normal linear system. However, similar to the
discussion in this subsection, for the normal linear system (\ref{normal}), we
can also set $F_{1}=0_{n\times n}$ such that the closed-loop system is
equivalent to an antilinear system. We may call this procedure as the
anti-linearization of the normal linear system (\ref{normal}). The
corresponding coupled matrix equations and their complete solutions can be
easily stated and are omitted for brevity.
\end{remark}

\section{\label{sec5}Sylvester and Stein Bimatrix Equations}

We will study in this section solutions to Sylvester and Stein bimatrix
equations. Since we are frequently dealing with the normal linear system
(\ref{normal}) and the antilinear system (\ref{antilinear}), for easy
reference, we make the following assumptions.

\begin{assumption}
\label{ass1}$A_{2}=0_{n\times n},F_{2}=0_{p\times p}$ and $C_{2}=0_{n\times
p},$ namely, $\{A_{1},A_{2}\}=\{A_{1},0_{n\times n}\},\{F_{1},F_{2}%
\}=\{F_{1},0_{p\times p}\}$ and $\{C_{1},C_{2}\}=\{C_{1},0_{n\times p}\}.$
\end{assumption}

\begin{assumption}
\label{ass2}$A_{1}=0_{n\times n},F_{1}=0_{p\times p}$ and $C_{1}=0_{n\times
p},$ namely, $\{A_{1},A_{2}\}=\{0_{n\times n},A_{2}\},\{F_{1},F_{2}%
\}=\{0_{p\times p},F_{2}\}$ and $\{C_{1},C_{2}\}=\{0_{n\times p},C_{2}\}.$
\end{assumption}

\subsection{\label{sec5.1}The Sylvester Bimatrix Equation}

In this subsection we discuss the Sylvester bimatrix equation
(\ref{bimatrixsyl}), which can also be written as coupled
matrix equations%
\begin{equation}
\left \{
\begin{array}
[c]{cl}%
C_{1} & =A_{1}X_{1}+A_{2}^{\#}X_{2}-(X_{1}F_{1}+X_{2}^{\#}F_{2}),\\
C_{2} & =A_{1}^{\#}X_{2}+A_{2}X_{1}-(X_{1}^{\#}F_{2}+X_{2}F_{1}).
\end{array}
\right.  \label{eqc1c2}%
\end{equation}

\begin{proposition}
\label{pp1}The Sylvester bimatrix equation (\ref{bimatrixsyl}) has a unique
solution if and only if
\begin{equation}
\lambda \left \{  A_{1},A_{2}\right \}  \cap \lambda \left \{  F_{1},F_{2}\right \}
=\varnothing. \label{eqeigjoint}%
\end{equation}
In this case, the unique solution is given by
\begin{equation}
\left \{  X_{1},X_{2}\right \}  =\left(  \sum \limits_{k=0}^{2p}\beta_{k}\left \{
A_{1},A_{2}\right \}  ^{k}\right)  ^{-1}\sum \limits_{k=1}^{2p}\beta_{k}\left \{
D_{1}(k),D_{2}(k)\right \}  , \label{equniquesylsolu}%
\end{equation}
where $\beta \left(  s\right)  =s^{2p}+\beta_{2p-1}s^{2p-1}+\cdots+\beta
_{1}s+\beta_{0}\ $is the characteristic polynomial of $\left \{  F_{1}%
,F_{2}\right \}  ,$ and, for $k\geq1,$%
\begin{equation}
\left \{  D_{1}(k),D_{2}(k)\right \}  =\sum \limits_{i=0}^{k-1}\left \{
A_{1},A_{2}\right \}  ^{i}\left \{  C_{1},C_{2}\right \}  \left \{  F_{1}%
,F_{2}\right \}  ^{k-1-i}. \label{eqd1kd2k}%
\end{equation}

\end{proposition}

Notice that the bimatrix series $\left \{  D_{1}(k),D_{2}(k)\right \}  $ in
(\ref{eqd1kd2k}) can also be defined in a recursive way:%
\begin{equation}
\left \{
\begin{array}
[c]{rl}%
\left \{  D_{1}(k+1),D_{2}(k+1)\right \}   & =\left \{  A_{1},A_{2}\right \}
\left \{  D_{1}(k),D_{2}(k)\right \}  +\left \{  C_{1},C_{2}\right \}  \left \{
F_{1},F_{2}\right \}  ^{k}\\
& =\left \{  A_{1},A_{2}\right \}  ^{k}\left \{  C_{1},C_{2}\right \}  +\left \{
D_{1}(k),D_{2}(k)\right \}  \left \{  F_{1},F_{2}\right \}  ,\;k\geq1,\\
\left \{  D_{1}(1),D_{2}(1)\right \}   & =\left \{  C_{1},C_{2}\right \}  .
\end{array}
\right.  \label{eq43}%
\end{equation}

\begin{corollary}
\label{ppsyl}Assume that%
\begin{equation}
\mu \left \{  A_{1},A_{2}\right \}  +\mu \left \{  F_{1},F_{2}\right \}  <0.
\label{eqmu1mu2}%
\end{equation}
Then the Sylvester bimatrix equation (\ref{bimatrixsyl}) has a unique solution
given by%
\begin{equation}
\left \{  X_{1},X_{2}\right \}  =\int_{0}^{\infty}\mathrm{e}^{t\left \{
A_{1},A_{2}\right \}  }\left \{  C_{1},C_{2}\right \}  \mathrm{e}^{t\left \{
F_{1},F_{2}\right \}  }\mathrm{d}t. \label{eq44}%
\end{equation}
Particularly, if the complex-valued linear system (\ref{sys}) with
$t\in \mathbf{R}^{+}$ is asymptotically stable, then the solution to the
Lyapunov bimatrix equation (\ref{bilya}) is given by
\begin{equation}
\left \{  P_{1},P_{2}\right \}  =\int_{0}^{\infty}\mathrm{e}^{t\left \{
A_{1},A_{2}\right \}  ^{\mathrm{H}}}\left \{  Q_{1},Q_{2}\right \}
\mathrm{e}^{t\left \{  A_{1},A_{2}\right \}  }\mathrm{d}t. \label{slya}%
\end{equation}

\end{corollary}

We now take a look at the following well-known Sylvester matrix equation%
\begin{equation}
A_{1}X-XF_{1}=C_{1},\label{syleq}%
\end{equation}
which was firstly studied by J. J. Sylvester \cite{sylvester84} and then by
several authors (for example, \cite{Hartwig72siam,jameson68siam}). The
following corollary reveals the relationship between the Sylvester bimatrix
equation (\ref{bimatrixsyl}) and the Sylvester matrix equation (\ref{syleq}).

\begin{corollary}
\label{normalsyl}The Sylvester matrix equation (\ref{syleq}) is solvable if
and only if the Sylvester bimatrix equation (\ref{bimatrixsyl}) under
Assumption \ref{ass1} is solvable. Particularly,

\begin{enumerate}
\item If $\left \{  X_{1},X_{2}\right \}  $ is a solution to (\ref{bimatrixsyl}%
), then $X=X_{1}$ is a solution to (\ref{syleq}), and if $X$ is a solution to
(\ref{syleq}), then $\{X_{1},X_{2}\}=\left \{  X,0_{n\times p}\right \}  $ is a
solution to (\ref{bimatrixsyl}).

\item If (\ref{bimatrixsyl}) has a unique solution $\left \{  X_{1}%
,X_{2}\right \}  ,$ then $X_{2}=0_{n\times p}$ and (\ref{syleq}) also has a
unique solution $X=X_{1}.$ If (\ref{syleq}) has a unique solution $X,$ then
(\ref{bimatrixsyl}) has a unique solution (given by $\left \{  X,0_{n\times
p}\right \}  $) if one of the following two conditions is satisfied%
\begin{align}
s  &  \in \lambda \left(  A_{1}\right)  \Longrightarrow s^{\#}\in \lambda \left(
A_{1}\right)  ,\label{eq27a}\\
s  &  \in \lambda \left(  F_{1}\right)  \Longrightarrow s^{\#}\in \lambda \left(
F_{1}\right)  . \label{eq27b}%
\end{align}

\end{enumerate}
\end{corollary}

Notice that equation (\ref{eq27a}) (equation (\ref{eq27b})) is satisfied if
$A_{1}$ ($F_{1}$) is a real matrix. We can check that, under Assumption
\ref{ass1}, the unique solution (\ref{equniquesylsolu}) to (\ref{bimatrixsyl})
coincides exactly with the unique solution to (\ref{syleq}) obtained in
\cite{jameson68siam}. We omit the details to save spaces.

\begin{remark}
\label{rm2}Under Assumption \ref{ass1}, it follows from $\lambda
\{A_{1},0_{n\times n}\}=\lambda(A_{1})\cup \lambda(A_{1}^{\#})$ that
(\ref{eqmu1mu2}) is equivalent to%
\begin{equation}
\mu \left(  A_{1}\right)  +\mu \left(  F_{1}\right)  <0. \label{eq64}%
\end{equation}
On the other hand, we have $\mathrm{e}^{t\left \{  A_{1},A_{2}\right \}  }=\{
\mathrm{e}^{A_{1}t},0_{n\times n}\}$ and $\mathrm{e}^{t\left \{  F_{1}%
,F_{2}\right \}  }=\{ \mathrm{e}^{F_{1}t},0_{p\times p}\}.$ Then applying
Corollary \ref{ppsyl} on (\ref{syleq}) gives the well-known result
\cite{Lancaster70siam}: if (\ref{eq64}) is satisfied, the Sylvester matrix
equation (\ref{syleq}) has a unique solution given by $X=\int_{0}^{\infty
}\mathrm{e}^{tA_{1}}C_{1}\mathrm{e}^{tF_{1}}\mathrm{d}t$.
\end{remark}

We next discuss the so-called conjugate-Sylvester matrix equation%
\begin{equation}
C_{2}=A_{2}^{\#}X-X^{\#}F_{2}, \label{sylconj}%
\end{equation}
which was firstly investigated in \cite{bhh87ctmt,bhh88siam} and then in
\cite{wz17book}.

\begin{corollary}
\label{conjsyl}The conjugate-Sylvester matrix equation (\ref{sylconj}) is
solvable if and only if the Sylvester bimatrix equation (\ref{bimatrixsyl})
under Assumption \ref{ass2}$\ $is solvable. Particularly,

\begin{enumerate}
\item If $\left \{  X_{1},X_{2}\right \}  $ is a solution to (\ref{bimatrixsyl}%
), then $X=X_{1}$ is a solution to (\ref{sylconj}), and if $X$ is a solution
to (\ref{sylconj}), then $\left \{  X_{1},X_{2}\right \}  =\left \{  X,0_{n\times
p}\right \}  $ is a solution to (\ref{bimatrixsyl}).

\item (\ref{bimatrixsyl}) has a unique solution (denoted by $\left \{
X_{1},X_{2}\right \}  $) if and only if (\ref{sylconj}) has a unique solution
(denoted by $X$). Moreover, $X_{2}=0_{n\times p}$ and $X=X_{1}$ with%
\begin{equation}
X_{1}=\gamma^{-1}\left(  A_{2}^{\#}A_{2}\right)  \left(  \sum \limits_{k=1}%
^{p}\gamma_{k}\sum \limits_{i=0}^{k-1}\left(  A_{2}^{\#}A_{2}\right)
^{i}\left(  C_{2}^{\#}F_{2}+A_{2}^{\#}C_{2}\right)  \left(  F_{2}^{\#}%
F_{2}\right)  ^{k-1-i}\right)  , \label{eqx2}%
\end{equation}
where $\gamma \left(  s\right)  =s^{p}+\gamma_{p-1}s^{p-1}+\cdots+\gamma
_{1}s+\gamma_{0}$ is a polynomial with real coefficients, defined by%
\begin{equation}
\gamma \left(  s\right)  =\left \vert sI_{p}-F_{2}^{\#}F_{2}\right \vert .
\label{eqg}%
\end{equation}

\end{enumerate}
\end{corollary}

We notice that (\ref{eqx2}) coincides with the solution obtained in
\cite{wz17book}. Moreover, Item 2 of Corollary \ref{conjsyl} seems better than
that of Corollary \ref{normalsyl} since it provides an \textquotedblleft if
and only if\textquotedblright \ condition.

\begin{remark}
\label{rm3}An explanation of the solution (\ref{eqx2}) is given as follows. We have
\begin{align*}
\left \{  D_{1}(2),D_{2}(2)\right \}   &  =\left \{  A_{2}^{\#}A_{2},0_{n\times
n}\right \}  \left \{  X_{1},X_{2}\right \}  -\left \{  X_{1},X_{2}\right \}
\left \{  F_{2}^{\#}F_{2},0_{p\times p}\right \} \\
&  =\left \{  A_{2}^{\#}A_{2}X_{1}-X_{1}F_{2}^{\#}F_{2},A_{2}A_{2}^{\#}%
X_{2}-X_{2}F_{2}^{\#}F_{2}\right \} \\
&  =\left \{  C_{2}^{\#}F_{2}+A_{2}^{\#}C_{2},0_{n\times p}\right \}.
\end{align*}
Therefore, if (\ref{bimatrixsyl}) has a
unique solution, then $X_{2}=0_{n\times p}$ and $X_{1}$ satisfies%
\begin{equation}
A_{2}^{\#}A_{2}X_{1}-X_{1}F_{2}^{\#}F_{2}=C_{2}^{\#}F_{2}+A_{2}^{\#}C_{2}.
\label{eq51}%
\end{equation}
The closed-form solution to (\ref{eq51}) is just (\ref{eqx2}) by using the
result in \cite{jameson68siam}.
\end{remark}

\subsection{The Stein Bimatrix Equation}

We now study the Stein bimatrix equation (\ref{bistein}), which is equivalent to the coupled matrix equations%
\begin{equation}
\left \{
\begin{array}
[c]{cc}%
X_{1} & =A_{1}X_{1}F_{1}+A_{2}^{\#}X_{2}F_{1}+A_{1}X_{2}^{\#}F_{2}+A_{2}%
^{\#}X_{1}^{\#}F_{2}+C_{1},\\
X_{2} & =A_{1}^{\#}X_{1}^{\#}F_{2}+A_{2}X_{2}^{\#}F_{2}+A_{1}^{\#}X_{2}%
F_{1}+A_{2}X_{1}F_{1}+C_{2}.
\end{array}
\right.  \label{steincoupled}%
\end{equation}

\begin{proposition}
\label{ppstein}The Stein bimatrix equation (\ref{bistein}) has a unique
solution if and only if
\begin{equation}
\lambda_{i}\left \{  A_{1},A_{2}\right \}  \lambda_{j}\left \{  F_{1}%
,F_{2}\right \}  \neq1,\forall i,j. \label{eq20}%
\end{equation}
In this case, the unique solution is given by
\begin{equation}
\left \{  X_{1},X_{2}\right \}  =\left(  \sum \limits_{k=0}^{2p}\beta_{k}\left \{
A_{1},A_{2}\right \}  ^{2p-i}\right)  ^{-1}\sum \limits_{k=1}^{2p}\beta
_{k}\left \{  A_{1},A_{2}\right \}  ^{2p-k}\left \{  D_{1}(k),D_{2}(k)\right \}  ,
\label{eq21}%
\end{equation}
where $\beta \left(  s\right)  =\beta^{2p}+\beta_{2p-1}s^{2p-1}+\cdots
+\beta_{1}s+\beta_{0}$ is the characteristic polynomial of $\left \{
F_{1},F_{2}\right \}  ,$ and, for $k\geq1,$%
\begin{equation}
\left \{  D_{1}(k),D_{2}(k)\right \}  =\sum \limits_{i=0}^{k-1}\left \{
A_{1},A_{2}\right \}  ^{i}\left \{  C_{1},C_{2}\right \}  \left \{  F_{1}%
,F_{2}\right \}  ^{i}. \label{eq22}%
\end{equation}

\end{proposition}

Notice that, similar to (\ref{eq43}), the bimatrix $\left \{  D_{1}%
(k),D_{2}(k)\right \}  $ in (\ref{eq22}) can also be defined in a recursive
way:%
\begin{equation}
\left \{
\begin{array}
[c]{rl}%
\left \{  D_{1}(k+1),D_{2}(k+1)\right \}   & =\left \{  A_{1},A_{2}\right \}
\left \{  D_{1}(k),D_{2}(k)\right \}  \left \{  F_{1},F_{2}\right \}  ,\;k\geq1,\\
\left \{  D_{1}(1),D_{2}(1)\right \}   & =\left \{  C_{1},C_{2}\right \}  .
\end{array}
\right.  \label{eqdd}%
\end{equation}
We can also state the following corollary that parallels to
corollary \ref{ppsyl}.

\begin{corollary}
\label{ppstein2}Assume that%
\begin{equation}
\rho \left \{  A_{1},A_{2}\right \}  \rho \left \{  F_{1},F_{2}\right \}  <1,
\label{eqrho}%
\end{equation}
and $\left \{  D_{1}(k),D_{2}(k)\right \}  $ is defined by (\ref{eq22}). Then
the Stein bimatrix equation (\ref{bistein}) has a unique solution given by%
\begin{equation}
\left \{  X_{1},X_{2}\right \}  =\sum \limits_{k=0}^{\infty}\left \{  A_{1}%
,A_{2}\right \}  ^{k}\left \{  C_{1},C_{2}\right \}  \left \{  F_{1}%
,F_{2}\right \}  ^{k}. \label{eq29}%
\end{equation}
Particularly, if the complex-valued linear system (\ref{sys}) with
$t\in \mathbf{Z}^{+}$ is asymptotically stable, then the solution to the
discrete-time Lyapunov bimatrix equation (\ref{bidilya}) is given by
\begin{equation}
\left \{  P_{1},P_{2}\right \}  =\sum \limits_{k=0}^{\infty}\left(  \left \{
A_{1},A_{2}\right \}  ^{\mathrm{H}}\right)  ^{k}\left \{  Q_{1},Q_{2}\right \}
\left \{  A_{1},A_{2}\right \}  ^{k}. \label{dlyasolu}%
\end{equation}

\end{corollary}

We now check the following well-known Stein matrix equation%
\begin{equation}
X=A_{1}XF_{1}+C_{1}, \label{stein}%
\end{equation}
which has been widely studied in the literature (see, for example,
\cite{jw03laa} and \cite{zld11laa}). The relationship between
(\ref{bistein}) and (\ref{stein}) can be made clear as follows.

\begin{corollary}
\label{normalstein}The Stein matrix equation (\ref{stein}) is solvable if and
only if the Stein bimatrix equation (\ref{bistein}) under Assumption
\ref{ass1} is solvable. Particularly,

\begin{enumerate}
\item If $\left \{  X_{1},X_{2}\right \}  $ is a solution to (\ref{bistein}),
then $X=X_{1}$ is a solution to (\ref{stein}), and if $X$ is a solution to
(\ref{stein}), then $\{X_{1},X_{2}\}=\left \{  X,0_{n\times p}\right \}  $ is a
solution to (\ref{bistein}).

\item If (\ref{bistein}) has a unique solution $\left \{  X_{1},X_{2}\right \}
,$ then $X_{2}=0_{n\times p}$ and (\ref{stein}) also has a unique solution
$X=X_{1}.$ If (\ref{stein}) has a unique solution $X,$ then (\ref{bistein})
has a unique solution (given by $\left \{  X,0_{n\times p}\right \}  $) if one
of the two conditions (\ref{eq27a})-(\ref{eq27b}) is satisfied.
\end{enumerate}
\end{corollary}

Under Assumption \ref{ass1}, one may check that the unique solution
(\ref{eq21}) to (\ref{bistein}) coincides with the unique solution to
(\ref{stein}) obtained in \cite{jw03laa}. The details are omitted.

\begin{remark}
Similar to Remark \ref{rm2}, under Assumption \ref{ass1}, we can show that
(\ref{eqrho}) is equivalent to%
\begin{equation}
\rho \left(  A_{1}\right)  \rho \left(  F_{1}\right)  <1. \label{eqrho1}%
\end{equation}
On the other hand, we have from (\ref{eqdd}) that $D_{2}(k)=0_{n\times p}$ and
$D_{1}(k+1)=A_{1}D_{1}\left(  k\right)  F_{1}$ with $D_{1}\left(  1\right)
=C_{1}.$ Then applying Corollary \ref{ppstein2} on (\ref{stein}) gives the
well-known result (see, for example, \cite{zld11laa}): if (\ref{eqrho1}) is
satisfied, the Stein matrix equation (\ref{stein}) has a unique solution given
by $X= \sum \nolimits_{k=0}^{\infty} A_{1}^{k}C_{1}F_{1}^{k}.$
\end{remark}

We next discuss the so-called conjugate-Stein matrix equation%
\begin{equation}
X=A_{2}X^{\#}F_{2}+C_{2}, \label{conj_stein}%
\end{equation}
which was firstly investigated in \cite{jw03laa} and was also studied in
several other papers, for example, \cite{wz17book} and \cite{zld11laa}, where
it was equivalently transformed into a normal Stein matrix equation in the
form of (\ref{stein}). The following corollary reveals the relationship
between (\ref{conj_stein}) and (\ref{bistein}).

\begin{corollary}
\label{conjstein}The conjugate-Stein matrix equation (\ref{conj_stein}) is
solvable if and only if the Stein bimatrix equation (\ref{bistein}) under
Assumption \ref{ass2} is solvable. Particularly,

\begin{enumerate}
\item If $\left \{  X_{1},X_{2}\right \}  $ is a solution to (\ref{bistein}),
then $X=X_{2}$ is a solution to (\ref{conj_stein}), and if $X$ is a solution
to (\ref{conj_stein}), then $\left \{  X_{1},X_{2}\right \}  =\left \{
0_{n\times p},X\right \}  $ is a solution to (\ref{bistein}).

\item (\ref{bistein}) has a unique solution (denoted by $\left \{  X_{1}%
,X_{2}\right \}  $) if and only if (\ref{conj_stein}) has a unique solution
(denoted by $X$). Moreover, $X_{1}=0_{n\times p}\ $and $X=X_{2}$ with%
\begin{equation}
X_{2}=\left(  \sum \limits_{k=0}^{p}\gamma_{k}(A_{2}A_{2}^{\#})^{p-k}\right)
^{-1}\left(  \sum \limits_{k=1}^{p}\gamma_{k}(A_{2}A_{2}^{\#})^{p-1}\left(
C_{2}+A_{2}C_{2}^{\#}F_{2}\right)  \left(  F_{2}^{\#}F_{2}\right)
^{k-1}\right)  . \label{eqx22}%
\end{equation}

\end{enumerate}
\end{corollary}

Similar to the situation in Subsection \ref{sec5.1}, Item 2 of Corollary
\ref{conjstein} is better than that of Corollary \ref{normalstein} since it
reveals an \textquotedblleft if and only if\textquotedblright \ relation.
Notice that (\ref{eqx22}) coincides with the solution obtained in
\cite{jw03laa}.

\begin{remark}
\label{rm5}Parallel to Remark \ref{rm3}, we give an explanation on
(\ref{eqx22}). We have
\begin{align*}
\left \{  X_{1},X_{2}\right \}   &  =\left \{  A_{1},A_{2}\right \}  ^{2}\left \{
X_{1},X_{2}\right \}  \left \{  F_{1},F_{2}\right \}  ^{2}+\left \{
D_{1}(2),D_{2}(2)\right \} \\
&  =\left \{  A_{2}^{\#}A_{2},0_{n\times n}\right \}  \left \{  X_{1}%
,X_{2}\right \}  \left \{  F_{2}^{\#}F_{2},0_{p\times p}\right \}  +\left \{
0_{n\times p},C_{2}+A_{2}C_{2}^{\#}F_{2}\right \} \\
&  =\left \{  A_{2}^{\#}A_{2}X_{1}F_{2}^{\#}F_{2},A_{2}A_{2}^{\#}X_{2}%
F_{2}^{\#}F_{2}+C_{2}+A_{2}C_{2}^{\#}F_{2}\right \}  .
\end{align*}
Therefore, if (\ref{bistein}) has a unique solution, we have $X_{1}=0_{n\times p}$ and $X_{2}$ satisfying%
\begin{equation}
X_{2}=A_{2}A_{2}^{\#}X_{2}F_{2}^{\#}F_{2}+C_{2}+A_{2}C_{2}^{\#}F_{2},
\label{eq54}%
\end{equation}
whose unique solution is exactly (\ref{eqx22}) by using the result in
\cite{jw03laa}. Moreover, (\ref{eqrho}) is equivalent to
\begin{equation}
\rho \left(  A_{2}A_{2}^{\#}\right)  \rho \left(  F_{2}^{\#}F_{2}\right)  <1.
\label{eq50}%
\end{equation}
Hence, under this condition, the unique solution to (\ref{eq54}) can also be
expressed by (see (72) in \cite{zld11laa}):%
\[
X_{2}=\sum \limits_{k=0}^{\infty}\left(  A_{2}A_{2}^{\#}\right)  ^{k}\left(
C_{2}+A_{2}C_{2}^{\#}F_{2}\right)  \left(  F_{2}^{\#}F_{2}\right)  ^{k}.
\]

\end{remark}

\section{\label{sec6}Example: Design of the Spacecraft Rendezvous System}

We use the spacecraft rendezvous system model to illustrate the effectiveness
of the proposed methods. The linearized equation of the spacecraft rendezvous
control system is known as the C-W equation \cite{cw60jas}%
\begin{equation}
\left \{
\begin{array}
[c]{l}%
\ddot{\xi}_{1}=2\omega \dot{\xi}_{2}+3\omega^{2}\xi_{1}+a_{1},\\
\ddot{\xi}_{2}=-2\omega \dot{\xi}_{1}+a_{2},\\
\ddot{\xi}_{3}=-\omega^{2}\xi_{3}+a_{3},
\end{array}
\right.  \label{cw}%
\end{equation}
where $\xi_{i},i=1,2,3,$ are relative positives of the chase spacecraft with
respect to the target, $a_{1},a_{2}$ and $a_{3}$ are the control accelerations
that the thrusts generate in the three directions, and $\omega$ is the orbit
rate of the target orbit, which is a known constant. For more information
about this model, see \cite{cw60jas} and the
references therein. Notice that (\ref{cw}) is a full-actuated system. To
demonstrate the design in a general case, we assume that $a_{1}=0$, since,
otherwise, the design is trivial as $a_{i}$ can be easily designed to cancel
all of the open-loop dynamics. Notice that in this case the system is still
controllable.

System (\ref{cw}) is exactly in the form of (\ref{second}) with $\xi=[\xi
_{1},\xi_{2},\xi_{3}]^{\mathrm{T}},$ $v=[a_{2},a_{3}]^{\mathrm{T}},n=3$, $q=2$
and $m=1$. Then, according to the development in Subsection \ref{sec3.3}, this
system can be written as the complex-valued linear system (\ref{sys}), where
$A_{i},B_{i},i=1,2,$ are computed according to (\ref{a12b12}) as%
\begin{align*}
A_{1}  &  =\left[
\begin{array}
[c]{ccc}%
\frac{3\omega^{2}\mathrm{j}}{2}-\frac{\mathrm{j}}{2} & \omega & 0\\
-\omega & -\frac{\mathrm{j}}{2} & 0\\
0 & 0 & -\frac{\omega^{2}\mathrm{j}}{2}-\frac{\mathrm{j}}{2}%
\end{array}
\right]  ,\;B_{1}=\left[
\begin{array}
[c]{c}%
0\\
\frac{\mathrm{j}}{2}\\
\frac{1}{2}%
\end{array}
\right]  ,\\
A_{2}  &  =\left[
\begin{array}
[c]{ccc}%
-\frac{3\omega^{2}\mathrm{j}}{2}-\frac{\mathrm{j}}{2} & -\omega & 0\\
\omega & -\frac{\mathrm{j}}{2} & 0\\
0 & 0 & \frac{\omega^{2}\mathrm{j}}{2}-\frac{\mathrm{j}}{2}%
\end{array}
\right]  ,\;B_{2}=\left[
\begin{array}
[c]{c}%
0\\
-\frac{\mathrm{j}}{2}\\
-\frac{1}{2}%
\end{array}
\right]  .
\end{align*}
The eigenvalue set of the open-loop system is known as $\{0,0,\pm
\omega \mathrm{j},\pm \omega \mathrm{j}\}$ \cite{cw60jas}. We are going to design
a full state feedback (\ref{eqfeedback}) such that the closed-loop system
possesses an eigenvalue set that is a shift of the open-loop system along the
real axis, say, $\mathit{\Gamma}=\{-\gamma,-\gamma,-\gamma \pm \omega
\mathrm{j},-\gamma \pm \omega \mathrm{j}\},$ where $\gamma>0$ is a real constant.
Thus we can choose%
\[
F=\mathrm{diag}\left \{  \left[
\begin{array}
[c]{cc}%
0 & 1\\
0 & 0
\end{array}
\right]  ,\left[
\begin{array}
[c]{cc}%
0 & \omega \\
-\omega & 0
\end{array}
\right]  ,\left[
\begin{array}
[c]{cc}%
0 & \omega \\
-\omega & 0
\end{array}
\right]  \right \}  -\gamma I_{6},
\]
and the unique bimatrix $\{F_{1},F_{2}\}$ satisfying (\ref{eqff}) can be
obtained as%
\[
F_{1}=\left[
\begin{array}
[c]{ccc}%
-\gamma & \frac{1}{2} & -\frac{\omega \mathrm{j}}{2}\\
0 & -\gamma & \frac{\omega}{2}\\
-\frac{\omega \mathrm{j}}{2} & -\frac{\omega}{2} & -\gamma
\end{array}
\right]  ,\;F_{2}=\left[
\begin{array}
[c]{ccc}%
0 & \frac{1}{2} & \frac{\omega \mathrm{j}}{2}\\
0 & 0 & -\frac{\omega}{2}\\
-\frac{\omega \mathrm{j}}{2} & \frac{\omega}{2} & 0
\end{array}
\right]  .
\]

The right-coprime bimatrix polynomials $\{N_{1}(s),N_{2}(s)\}$ and
$\{D_{1}(s),D_{2}(s)\}$ satisfying (\ref{coprime}) can be computed as%
\[
\left \{
\begin{array}
[c]{rl}%
N_{1}\left(  s\right)  & =\left[
\begin{array}
[c]{c}%
\omega s\left(  1+\mathrm{j}s\right) \\
-\frac{\left(  1+\mathrm{j}s\right)  \left(  3\omega^{2}-s^{2}\right)  }{2}\\
\frac{s}{2}-\frac{\mathrm{j}}{2}%
\end{array}
\right]  \triangleq \sum \limits_{i=0}^{3}N_{1i}s^{i},\\
N_{2}\left(  s\right)  & =\left[
\begin{array}
[c]{c}%
-\omega s\left(  -1+\mathrm{j}s\right) \\
\frac{\mathrm{j}\left(  3\omega^{2}-s^{2}\right)  \left(  s+\mathrm{j}\right)
}{2}\\
-\frac{s}{2}-\frac{\mathrm{j}}{2}%
\end{array}
\right]  \triangleq \sum \limits_{i=0}^{3}N_{2i}s^{i},
\end{array}
\right.
\]
and%
\[
\left \{
\begin{array}
[c]{rl}%
D_{1}\left(  s\right)  & =\frac{1}{2}\left(  s^{2}+1\right)  \left(
\omega^{2}+s^{2}\right)  \triangleq \sum \limits_{i=0}^{4}D_{1i}s^{i},\\
D_{2}\left(  s\right)  & =\frac{1}{2}\left(  s^{2}-1\right)  \left(
\omega^{2}+s^{2}\right)  \triangleq \sum \limits_{i=0}^{4}D_{2i}s^{i}.
\end{array}
\right.
\]
Then, according to Theorem \ref{th1}, complete solutions to the associated
generalized Sylvester bimatrix equation (\ref{syl}) are given by%
\[
\left \{
\begin{array}
[c]{rl}%
\left \{  X_{1},X_{2}\right \}  & =\sum \limits_{i=0}^{3}\left \{  N_{1i}%
,N_{2i}\right \}  \left \{  Z_{1},Z_{2}\right \}  \left \{  F_{1},F_{2}\right \}
^{i},\\
\left \{  Y_{1},Y_{2}\right \}  & =\sum \limits_{i=0}^{4}\left \{  D_{1i}%
,D_{2i}\right \}  \left \{  Z_{1},Z_{2}\right \}  \left \{  F_{1},F_{2}\right \}
^{i},
\end{array}
\right.
\]
where $Z_{i}\in \mathbf{C}^{1\times3},i=1,2,$ are arbitrary matrices.
Particularly, if we choose%
\[
Z_{1}=\left[
\begin{array}
[c]{ccc}%
1+\mathrm{j} & 0 & 0
\end{array}
\right]  ,\;Z_{2}=\left[
\begin{array}
[c]{ccc}%
0 & 0 & 1
\end{array}
\right]  ,
\]
the state feedback gain bimatrix $\{K_{1},K_{2}\}$ can be obtained according
to (\ref{eqsylk}) as (we omit to display the variables $\{X_{1},X_{2}\}$ and
$\{Y_{1},Y_{2}\}$)%
\[
\left \{
\begin{array}
[c]{rl}%
K_{1} & =\left[
\begin{array}
[c]{ccc}%
\frac{k_{11}}{12\omega^{3}} & \frac{k_{12}}{6\omega^{2}} & -\frac
{\gamma \left(  2+\gamma \mathrm{j}\right)  }{2}%
\end{array}
\right]  ,\\
K_{2} & =\left[
\begin{array}
[c]{ccc}%
-\frac{k_{21}}{12\omega^{3}} & \frac{k_{22}}{6\omega^{2}} & \frac
{\gamma \left(  2+\gamma \mathrm{j}\right)  }{2}%
\end{array}
\right]  ,
\end{array}
\right.
\]
where
\[
\left \{
\begin{array}
[c]{l}%
k_{11}=\gamma^{4}\mathrm{j}-12\gamma^{3}\omega^{2}+19\gamma^{2}\omega
^{2}\mathrm{j}+\gamma^{2}-42\gamma \omega^{4}+6\gamma \omega^{2}\mathrm{j}%
+4\omega^{2},\\
k_{21}=-\gamma^{4}\mathrm{j}+12\gamma^{3}\omega^{2}-19\gamma^{2}\omega
^{2}\mathrm{j}+\gamma^{2}+42\gamma \omega^{4}+6\gamma \omega^{2}\mathrm{j}%
+4\omega^{2},\\
k_{12}=\gamma^{4}+\gamma^{2}\omega^{2}-\gamma^{2}\mathrm{j}+12\gamma \omega
^{2}\mathrm{j}-\omega^{2}\mathrm{j,}\\
k_{22}=\gamma^{4}+\gamma^{2}\omega^{2}+\gamma^{2}\mathrm{j}+12\gamma \omega
^{2}\mathrm{j}+\omega^{2}\mathrm{j.}%
\end{array}
\right.
\]
Finally, the resulting controller can be implemented according to (\ref{eqv}),
which is physically realizable.

\section{\label{sec7}Conclusion}

This paper has studied several kinds of linear bimatrix equations, whose
coefficients are bimatrices that were introduced in our early studies. These
equations arise from full state feedback pole assignment and stability
analysis of complex-valued linear systems. Explicit solutions to these linear
bimatrix equations are established. Particularly, explicit solutions are
provided for the case that the coefficients of the bimatrix equations are
determined by the so-called antilinear systems. Explicit solutions are then
used to solve the pole assignment problem for complex-valued linear systems,
particularly, for second-order linear systems that can be easily converted
into complex-valued linear systems. The spacecraft rendezvous control system
is then used to demonstrate the obtained theoretical results. The results in
this paper can be readily extended to linear bimatrix equations associated
with complex-valued descriptor linear systems and high-order complex-valued
linear systems.

\end{document}